\documentclass[11pt]{article}
\usepackage[margin=2.95cm]{geometry}
\usepackage{amssymb}
\usepackage{latexsym}
\begin{document}

\textheight=198mm
\textwidth=130mm
\hyphenation{group-oid group-oids quasi-hori-zon-tal e-di-tion}
\newtheorem{thm}{Theorem}[section]
\newtheorem{thmdef}[thm]{Theorem/Definition}
\newtheorem{prop}[thm]{Proposition} \newtheorem{lemma}[thm]{Lemma}
\newtheorem{cor}[thm]{Corollary} \newtheorem{dfn}[thm]{Definition}
\newtheorem{conj}[thm]{Conjecture}
\newtheorem{axiom}[thm]{Axiom} \newtheorem{rmk}[thm]{Remark}
\newtheorem{ex}[thm]{Example} \newtheorem{question}[thm]{Question}
\newtheorem{problem}[thm]{Problem}
\newcommand {\pf}{\noindent{\bf Proof.}\ }
\newcommand{\complex}{{\mathbb C}}
\newcommand{\naturals}{{\mathbb N}}
\newcommand{\reals}{{\mathbb R}}
\newcommand{\torus}{{\mathbb T}}
\newcommand{\integers}{{\mathbb Z}}
\newcommand{\trans}{{\pitchfork}}
\newcommand{\sgn}{{\rm sgn~}}
\newcommand{\nor}{{\rm nor}}
\newcommand{\ver}{{\rm ver}}
\newcommand{\hor}{{\rm hor}}
\newcommand{\Hom}{{\rm Hom}}
\newcommand{\vol}{{\rm vol}}
\newcommand{\cok}{{\rm coker}}
\newcommand{\Tr}{{\rm Tr}}
\newcommand{\covder}{{\bf D}}
\newcommand{\grad}{\nabla}
\newcommand{\ave}{\cala}
\newcommand{\doub}{\calv}
\newcommand{\gammabar}{{\overline{\Gamma}}}
\newcommand{\ebar}{{\overline{E}}}
\newcommand{\fbar}{{\overline{F}}}
\newcommand{\mbar}{{\overline{M}}}
\newcommand{\nbar}{{\overline{N}}}
\newcommand{\pbar}{{\overline{P}}}
\newcommand{\tbar}{{\overline{T}}}
\newcommand{\xbar}{{\overline{X}}}
\newcommand{\ybar}{{\overline{Y}}}
\newcommand{\zbar}{{\overline{Z}}}
\newcommand{\proj}{{\mathbb P}}	
\newcommand{\cala}{{\cal A}}
\newcommand{\calb}{{\cal B}}
\newcommand{\calc}{{\cal C}}
\newcommand{\cald}{{\cal D}}
\newcommand{\cale}{{\cal E}}
\newcommand{\calf}{{\cal F}}
\newcommand{\calg}{{\cal G}}
\newcommand{\calh}{{\cal H}}
\newcommand{\calj}{{\cal J}}
\newcommand{\call}{{\cal L}}
\newcommand{\calm}{{\cal M}}
\newcommand{\calo}{{\cal O}}
\newcommand{\calp}{{\cal P}}
\newcommand{\cals}{{\cal S}}
\newcommand{\calt}{{\cal T}}
\newcommand{\calv}{{\cal V}}
\newcommand{\calu}{{\cal U}}
\newcommand{\calw}{{\cal W}}
\newcommand{\boldf}{{\mathbf f}}
\newcommand{\boldg}{{\mathbf g}}
\newcommand{\del}{\partial}
\newcommand{\qed}{\begin{flushright} $\Box$\ \ \ \ \ \end{flushright}}
\newcommand{\lqed}{\begin{flushright} $\triangle$\ \ \ \ \ \end{flushright}}
\newcommand{\half}{\textstyle{\frac{1}{2}}}
\newcommand{\frakb}{\mathfrak{b}}
\newcommand{\frakg}{\mathfrak{g}}
\newcommand{\frakh}{\mathfrak{h}}
\newcommand{\frakk}{\mathfrak{k}}
\newcommand{\frakm}{\mathfrak{m}}
\newcommand{\frakt}{\mathfrak{t}}
\newcommand{\boe}{\mathbf{e}}
\newcommand{\bof}{\mathbf{f}}
\newcommand{\bog}{\mathbf{g}}
\newcommand{\rel}{\mathbf{REL}}
\newcommand{\linsymp}{\mathbf{LINSYMP}}
\newcommand{\symp}{\mathbf{SYMP}}
\newcommand{\four}{\mathbf{FOUR}}
\newcommand{\micsymp}{\mathbf{MICSYMP}}
\newcommand{\man}{\mathbf{MAN}}
\newcommand{\rat}{\mathbf{RAT}}
\newcommand{\RAT}{\mathbf{RAT}}
\newcommand{\mic}{\mathbf{MIC}}
\newcommand{\donsharp}{{\rm Don}^\sharp}
\newcommand{\arrows}{\,\lower1pt\hbox{$\longrightarrow$}\hskip-.24in\raise2pt
             \hbox{$\longrightarrow$}\,}
\newcommand{\defequal}{\stackrel{\mbox {\tiny {def}}}{=}}
\newcommand{\Lag}{\mathrm{Lag}}
\newcommand{\tr}{\mathrm{t}}
\newcommand{\red}{\mathrm{red}}
\newcommand{\topp}{\mathrm{top}}
\newcommand{\bDelta}{\mathbf{\Delta}}
\newcommand{\bepsilon}{\mathbf{\epsilon}}
\newcommand{\wedgetop}{\bigwedge^{\rm top}}

\title{{\bf Symplectic categories}}
\author
{Alan
Weinstein
\thanks{Research partially supported by NSF Grant
DMS-0707137
\newline \mbox{~~~~}MSC2010 Subject Classification Number: 
53D12 (Primary), 81S10 (Secondary).
\newline \mbox{~~~~}Keywords: symplectic manifold, lagrangian
submanifold, canonical relation, category, quantization}
\\Department of Mathematics\\ University of California\\Berkeley, CA
94720 USA\\ {\small(alanw@math.berkeley.edu)}}
\maketitle

\begin{abstract} 
Quantization problems suggest that the category of symplectic manifolds and
symplectomorphisms be augmented by the inclusion of canonical
relations as morphisms.
These relations compose well when a
transversality condition is 
satisfied, but the failure of the most general
 compositions to be smooth manifolds
means that the canonical relations do not comprise the
morphisms of a category. 

We discuss several existing and potential
remedies to the nontransversality problem.  Some of these involve
restriction to classes of lagrangian submanifolds for which the
transversality property automatically holds.  Others involve allowing
lagrangian ``objects'' more general than submanifolds.


\end{abstract}

\section{Introduction}
\label{sec-intro}
This paper is based on  two lectures given at the Geometry
Summer School at the Istituto Superior Tecnico in Lisbon, in July of 2009.
We describe several ongoing efforts to build
categories whose 
objects are symplectic
manifolds and whose morphisms are canonical relations.  
The lectures also included a discussion of categories whose
hom-objects are symplectic manifolds, with the composition
of morphisms a canonical relation, but we do not include that topic in
this paper.

The material presented here is based on the work of many people,
including the author's work in progress with
Alberto Cattaneo, Benoit Dherin,
Shamgar Gurevich, and Ronny Hadani.

\subsection{Canonical relations as morphisms in a category}

A {\bf canonical relation} between symplectic manifolds $M$ and $N$
is, by definition, a lagrangian submanifold of $M \times \overline{N}$, where 
$\overline{N}$ is $N$ with its symplectic structure multiplied by
$-1$.  For example, the graph of a symplectomorphism is a canonical
relation, as is any product of lagrangian submanifolds in $M$ and $N$.

In his work on Fourier integral operators, H\"ormander
\cite{ho:fourier}, following Maslov \cite{ma:theory}, 
observed that, under a transversality assumption, the set-theoretic composition 
of two canonical relations is again a
canonical relation, and that this composition is a 
``classical limit'' of the composition of certain operators.  

Shortly thereafter, Sniatycki and Tulczyjew \cite{sn-tu:generating}
defined {\bf symplectic
relations} as isotropic\footnote{The calculus of coisotropic relations
does not seem to have been introduced until much later, in
\cite{we:coisotropic}.}
submanifolds of products and showed that this
class of relations was closed under ``clean'' composition (see Section
\ref{sec-relations} 
below).  They also observed that the natural relation between a
symplectic manifold and the quotient of a submanifold by the kernel
of the pulled-back symplectic form is a symplectic relation.

Following in part some (unpublished) ideas of the author, 
Guillemin and Sternberg \cite{gu-st:problems} observed that the linear
canonical relations (i.e., lagrangian subspaces of products of
symplectic vector spaces) could be considered as the morphisms of a
category, and they constructed a partial quantization of this
category (in which lagrangian subspaces are enhanced by
half-densities.  
The automorphism groups in this category are the linear
symplectic groups, and the restriction of the 
Guillemin-Sternberg  quantization to each such group is a metaplectic
representation.   On the other hand, the quantization of certain
compositions of canonical relations leads to ill-defined operations at
the quantum level, such as the evaluation of a delta ``function'' at
its singular point, or the multiplication of delta functions.

The quantization of the linear symplectic category was part of a
larger project of 
quantizing canonical relations (enhanced with extra structure, such as
half-densities) in a functorial way, and this program
was set out more formally by the present author in  
\cite{we:symplectic geometry} and \cite{we:symplectic category}.  It
was advocated there that canonical relations should be
considered as the morphisms of a ``category'', and that quantization
should be a functor from there to a category of linear spaces
and linear maps, consistent with some additional structures.
 The word ``category'' appears in quotation marks above
because the composition of canonical relations can fail to be a
canonical relation, as will be explained in detail below, so we do not
have a category.  Briefly,
there are two problems.

\begin{itemize}
\item The composition of two canonical relations may fail even to be a
  manifold.
\item In the linear symplectic category, where each space of morphisms
  has a natural topology as a lagrangian grassmannian manifold, the
  composition operation is discontinuous.
\end{itemize}

These two problems are both related to the possible failure of 
 a transversality condition.   It
is in general hard to 
 remedy this by imposing  conditions on
{\em individual} canonical relations which are significantly weaker than
local invertibility, though
we do present in Section \ref{sec-micromorphisms} 
below such a condition 
which yields a 
category whose objects are {\em germs} of symplectic manifolds around
lagrangian submanifolds.   Otherwise, we must do something more
drastic, using objects more general than canonical relations as the
morphisms in our category, or modifying the notion of category itself.
  
We begin this paper with a general discussion of the category in which
the morphisms are arbitrary relations between sets.  We then take a
first look at the composition of linear and nonlinear canonical
relations and at
 reduction by coisotropic submanifolds of symplectic manifolds.
Next, we present a simple idea of Wehrheim and Woodward
\cite{we-wo:functoriality}, who embed the
canonical relations in a true category $\symp$ in what is in some sense the
optimal way.  Then we describe the category of microfolds
\cite{ca-dh-we:symplectic}, 
whose
objects are germs of symplectic manifolds around lagrangian
submanifolds.  At this point, we are already outside the setting in
which the morphisms are maps between sets, but we go even further in
the next section by describing the construction by Wehrheim and
Woodward \cite{we-wo:functoriality}
of a 2-category of which $\symp$ could be thought of as the
``coarse moduli space''.  After that, we review and extend an idea of
Sabot \cite{sa:electrical}, who
deals with the discontinuities of the composition of linear canonical
relations by forming the closure of the graph of the composition
operation.  The result is a multiple-valued operation whose graph is
an algebraic subvariety of a product of Grassmann manifolds.
Finally, we attempt to unify and extend the examples above by using the
language of simplicial spaces, in which categories and groupoids
appear as objects satisfying special ``Kan conditions'' (see, for
example, \cite{zh:ngroupoids}).
For each of the remedies above, we briefly
discuss the quantization
of the resulting structure. 



\bigskip
\noindent
{\bf Acknowledgments.}  I would like to thank my hosts in the group,
Analyse Alg\'ebrique, at the Institut Math\'ematique de
Jussieu (Paris), where this paper was written.  
I would also like to thank 
John Baez, Christian Blohmann, Ralph Cohen, Benoit Dherin,
Dan Freed, Dmitry Roytenberg, Graeme Segal, 
Katrin Wehrheim, Chris Woodward, and Chenchang Zhu for helpful suggestions.

\section{Relations and their composition}
\label{sec-relations}
We begin with the category $\rel$ whose objects are sets
and for which the morphism space $\rel(X,Y)$ is simply the
set of all subsets of $X\times Y$.  We consider each
such relation $f$ as a morphism \emph{to} $X$ \emph{from} $Y$.  Linear
[affine] subspaces of vector [affine] spaces form a subcategory of
$\rel$.  

The natural exchange mappings $X\times Y \rightarrow Y\times X$ define a
contravariant {\bf
  transposition} functor $f\mapsto f^\tr$
 from $\rel$ to itself.

For any relation to $X$ from $Y$,  $X$  is the {\bf target} and $Y$
the {\bf source}.  The image $f(Y)$ of $f$ under projection to $X$ is the
{\bf range} of $f$, and the image  $f^\tr(X) \subseteq Y$ is the {\bf
  domain} of $f$.
$f$ is {\bf surjective} if its range equals its target, and {\bf
  cosurjective} if its domain equals its source (i.e. if it is
``defined everywhere'').  For any  $y \in Y$, $f(y)$  denotes
the image of $f$ on $\{y\}$, i.e. the subset $\{x\in X|(x,y)\in f\}$ of $X$.

The composition $f\circ g$ of $f\in \rel(X,Y)$ with $g\in \rel(Y,Z)$
is $$\{(x,z)|\exists y\in Y ~{\rm such~that} (x,y)\in f ~{\rm and}~
(y,z) \in g\}.$$ 
It is useful to think of this as the result of a sequence
of three operations: first, form the product $f \times g \subset X\times
Y\times Y\times Z$; second, intersect it with $X \times
\Delta_Y \times Z$, where $\Delta_Y$ is the diagonal in $Y\times Y$, to
obtain the fibre product $f\times_Y g$;
third, project this intersection into $X\times Z$.

When $X$
and $Y$ are manifolds and $f$ is a (locally closed) 
submanifold of
$X\times Y$,  $f$ is a {\bf smooth relation}.
When $f\in\rel(X,Y)$ and $g\in\rel(Y,Z)$ are smooth,  the 
pair $(f,g)$ is {\bf transversal} if 
$f\times g$ is
transversal to $X \times
\Delta_Y \times Z$, so that their intersection $f\times_Y g$ is again a
manifold.   
A transversal pair is {\bf strongly transversal}, and we will
write $f \trans g$, if
the projection map from $f\times_Y g$ to $X\times Z$ is an embedding
onto a locally closed submanifold, 
in which case the image $f\circ g$ is again a smooth relation.   When a
pair is not strongly transversal, its composition may fail to be a
submanifold, so the smooth relations do not form a
subcategory of $\rel$.   

A pair $(f,g)$ of linear or affine relations is 
transversal if and only if the domain of $f$ is transversal to the
range of $g$ as subspaces of $Y$, in which case the pair is
necessarily strongly transversal.  In particular, $f\trans g$
 whenever $f$ is cosurjective or $g$ is surjective.
Transversality of smooth $(f,g)$ is
detected by the same criterion, applied fibrewise to the
tangent relations $Tf$ and $Tg$.  If $f$ is
 the graph of a smooth mapping to $X$ from $Y$, then $f\trans g$
is a transversal pair for any smooth relation $g$ to $Y$ from $Z$.  
Similarly, if $g$ is the transpose of the
graph of a smooth mapping from $Y$
to $Z$, then $f\trans g$ for every
smooth relation $f$ to $X$ from $Y$.    In particular, the category
of smooth manifolds and smooth mappings is a subcategory of $\rel$.  

There is a condition weaker than transversality which, together with
an embedding condition, still insures that the composition of two
smooth relations is again smooth.  The pair $(f,g)$ is {\bf clean} 
if the fibre product $f\times_Y g$ is a submanifold of $X\times
Y\times Y\times Z$, and if the natural inclusion of $T(f\times_Y g)$ in
the fibre product tangent bundle $Tf\times_{TY} Tg$ is an equality
(equivalently, if $T(f\circ g) = Tf \circ Tg$),
and if the differential of the 
projection from $f\times_Y g$ to $X\times Z$ has constant
rank.   If this projection is a submersion onto
a locally closed submanifold of
 $X\times Z$, then $f\circ g$ 
is again a smooth relation, and the pair $(f,g)$ is {\bf
  immaculate}.   For example, any composable pair of linear or affine relations
is immaculate.

When $X$, $Y$, and $Z$ are symplectic manifolds, 
$C=X \times \Delta_Y \times Z$ is a coisotropic submanifold of 
$X \times \overline{Y} \times Y \times \overline{Z}$, and the leaves
of the characteristic foliation are connected components of the fibres of the
projection from $C$ to $X \times \overline{Z}$.  It follows that, for
canonical relations,  
the constant rank condition in the definition of a clean pair
follows from the other conditions.  For any transversal pair,
the projection from $f\times_Y g$ to $X \times \overline
Z$ is an immersion.

Finally, we note that any (lagrangian) submanifold $L$ of a
(symplectic) manifold $X$ may be thought of as a smooth (canonical)
relation to $X$ from a point or to a point from $X$.  Although points
are 
neither initial
nor terminal objects in our categories of relations, they still
play a special role.  For instance, the composition of $L$ to a point
from $X$ with $L'$ to $X$ from a point is nonempty if and only if $L$
is disjoint from $L'$.  Upon quantization, a point usually becomes the
scalars $\complex$, and the intersection of two lagrangian submanifolds
represents geometrically the inner product of quantum states to which they
correspond.  

\section{Canonical relations and coisotropic reduction}
\label{sec-reduction}
There are several connections between canonical relations and the
reduction of symplectic manifolds by coisotropic
submanifolds.   

\subsection{The linear case}
We begin with the linear case, though much of what we write here
applies immediately to the case of manifolds.  (See the following
section for further details.)

If $C$ is a coisotropic subspace of a symplectic vector space $X$, the
quotient $X_C$ of $C$ by the kernel $C^\perp$ of the induced bilinear
form carries a natural symplectic structure and is called the {\bf
  reduced space}.  It is connected to $X$ by 
the canonical relation 
$$r_c = \{(x,y)\in X_C \times X | y \in C ~{\rm and}~ x = [y]\},$$ where $[y]$ is the
equivalence class of $y$ modulo $C^\perp$.  The composition of $r_C$ with a
lagrangian subspace $L$ in $X$ (a linear canonical relation to $X$
from a point) gives the reduced lagrangian subspace $L_C = (L\cap C)/(L\cap
C^\perp)$ in $X_C$.  
The composition is transversal, hence strongly transversal, exactly when $L$
is transversal to $C$.  We will refer to ``transversal reduction'' in
this situation.  The  operation $\Lag (X_C)\leftarrow \Lag(X)$ given
by composition with $r_C$
will be denoted by $R_C$.

For any linear canonical relation $f\in \rel(X,Y)$, the range and domain 
are
coisotropic subspaces of $X$ and $Y$ respectively, and $f$ induces an
isomorphism $  X_{f^\tr (Y)} \stackrel{f_\red}{\longleftarrow} Y_{f(X)}$ between
the reduced spaces, giving a natural factorization
$$f = r_{f(Y)}^\tr \circ f_\red \circ r_{f^\tr(X)} $$ of any linear
canonical relation as the (transversal!) composition 
of a transposed reduction, a
symplectomorphism, and a reduction.  In particular, any surjective
linear canonical relation is essentially a reduction, and any
cosurjective one is
essentially a transposed reduction.

The composition of linear canonical relations is itself an instance of
reduction  As we have already mentioned in the previous section,
$X \times \Delta_Y \times Z$ is a coisotropic subspace
of $X \times \overline{Y} \times Y \times \overline{Z}$,
and $$(X\times \Delta_Y \times Z)^\perp 
= \{0_X\} \times \Delta_Y \times \{0_Y\}, $$ 
so
$(X \times \overline{Y} \times Y \times \overline{Z})_{X \times
  \Delta_Y \times Z}$ is naturally isomorphic to $X\times
\overline{Z}$.   Under this isomorphism, the composed relation
is the reduction of the product $f\times g$, and the composition
is transversal if and only if the reduction is.

\subsection{The nonlinear case}

We turn now to the case of manifolds.
Any coisotropic submanifold $C$ of a symplectic manifold $X$ carries
a {\bf characteristic distribution} $TC^\perp \subseteq TC$ 
which is, by definition, the kernel of the 
pullback to $C$ of the symplectic form on $X$.
$TC^\perp$ consists  of the values of hamiltonian vector
fields whose hamiltonians vanish on $C$ 
 and is always an
integrable distribution, hence tangent to a foliation which we will
denote by $\calc ^\perp$.  If $\calc ^\perp$ is {\bf simple} in the
sense that its leaves are the fibres of a submersion, then the leaf
space $X_C = C/\calc^\perp$ is again a symplectic manifold, and 
the reduction operation 
$$r_C = 
\{(x,y)\in X_C \times X | x = [y], y \in C \},$$ where $[y]$ is the
leaf of $\calc ^\perp$ containing $y$, is a canonical relation to
$X_C$ from $X$.   

Let $\pi_C$ be the projection from $C$ to $X_C$.  For any lagrangian
$X\subset C$ (a canonical relation to $X$ from a point), the
composition $r_C\circ L = r_C(L)$ is just the reduction  $L_C = \pi_C(L\cap C)$.  
The composition is transversal when $L$ is transversal to $C$, in
which case the restriction of $\pi_C$ to $L\cap C$ is an immersion
onto $L_C$,
but it is not necessarily injective.  Furthermore, if $L$ is not
transversal to $C$, then $L\cap C$ may not even be a manifold.  There, 
a reduction operation $R_C$ to lagrangian submanifolds of $X_C$ from
the lagrangian submanifolds of $X$ is not just discontinuous as in the
linear case, but not even defined everywhere.   



\section{The Wehrheim-Woodward category}
\label{sec-ww1}

Wehrheim and Woodward \cite{we-wo:functoriality} begin with the 
following construction to circumvent the problem of bad compositions.
(Our terminology and notation throughout this section differ somewhat
from theirs.)   The result is in some sense the minimal way to produce
a category whose morphisms include all of the canonical relations.  

\begin{dfn} The Wehrheim-Woodward category $\symp$ is the category
  whose objects are symplectic manifolds and whose morphisms are
  generated by the canonical relations, subject to the relation that
  the composition of $f$ and $g$ in $\symp$ is equal to the
  composition in $\rel$ when $f\trans g$.
\end{dfn}

More explicitly, as in \cite{we-wo:functoriality}, we may
begin with the category whose
objects are symplectic manifolds and whose morphisms are sequences
$(f_1,\ldots,f_r)$ of canonical relations which are  composable
in $\rel$.  We also include an empty sequence for each object, which
functions as an identity morphism.  
Composition is given by concatenation of
sequences.   Set-theoretic composition of relations defines a functor
from this category to $\rel$.

Now introduce the smallest equivalence relation which is closed
under composition from both sides, for which $(f,g)$ is equivalent
to $fg$ when $(f,g)$ is a strongly transversal pair, and for which each
empty sequence is equivalent to the graph of the identity map on the
corresponding object.  The equivalence classes are the morphisms in
$\symp$, and the composition functor above descends to give a functor
from $\symp$ to $\rel$.  It follows that distinct canonical relations,
considered sequences with a single entry, give distinct morphisms in $\rel$.
The identity morphisms are the (equivalence classes of the) diagonals
$\Delta_Y\subset Y\times \overline{Y}$.  A morphism in $\symp$ is
called a {\bf generalized lagrangian correspondence} in
\cite{we-wo:functoriality}, and a generalized lagrangian
correspondence to $X$ from a point is a {\bf generalized lagrangian
submanifold}.  

$\symp$ is characterized by the universal property that 
any map $\Phi$ to any category $\calc$ from 
symplectic manifolds and canonical relations such that
$\Phi(f\circ g) = \Phi(f)\Phi(g)$ whenever $f\trans g$
 factors uniquely through $\symp$.    It is thus tempting to
think of $\symp$ as a universal quantization category.  On the other
hand,  quantization of canonical relations by operators on function
spaces
requires
enhancement of the morphisms by some extra structure, such as
half-densities or half-forms.  Thus, it is natural to try to extend
$\symp$ and its variants by building larger categories with forgetful
functors to $\symp$.   It will also be important to extend to such
categories the basic operations on canonical relations, such as
transpose and cartesian products.

The construction of $\symp$ is merely the beginning of  
what
Wehrheim and Woodward do in \cite{we-wo:functoriality}.  Imposing
topological conditions (involving Chern classes, Maslov classes, etc.)
on the symplectic manifolds and canonical relations they define a
certain subcategory\footnote{Actually, there are two categories, one
  for exact symplectic manifolds and one for monotone ones.}  $\symp '$
of $\symp$ whose objects are ``admissible'' symplectic manifolds and
whose morphisms are generated by ``admissible'' canonical relations.  
For each admissible manifold $X$, they
construct a Donaldson-Fukaya category $\donsharp(X)$ whose objects are
generalized lagrangian submanifolds in $X$ and whose morphism spaces
are Floer cohomology groups.  Composition of morphisms 
involves counting pseudoholomorphic curves.
For each admissible relation
$f$ to $X$ from $Y$, they construct a functor $\donsharp(f)$ to 
$\donsharp(X)$ from  $\donsharp (Y) $ such that,
when $(f,g)$ is a strongly
transversal pair,
the functors $\donsharp(f\circ g)$ and $\donsharp (f) \circ
\donsharp(g)$ are naturally equivalent.  
$\donsharp$ then extends to a functor from 
$\symp'$ to the category whose objects are categories and whose
morphisms are natural equivalence classes of functors.   Since the
target category of this functor has an additive structure, we may view
$\donsharp$ as a kind of ``quantization'' of $\symp'$.  

We will continue our discussion of $\symp'$ and its quantization in
Section \ref{sec-ww2}.

\section{Cotangent lifts and symplectic micromorphisms}
\label{sec-micromorphisms}
There is a natural but quite limited collection of symplectic
manifolds and canonical relations which form a subcategory of
$\symp$.  It is the image of a contravariant functor $T^*$ from the category 
$\man$ of smooth manifolds and smooth maps.  
Namely, for every smooth manifold $M$, $T^*M$ is
its cotangent bundle with the canonical symplectic structure, and for
every smooth map\footnote{Note that
 the arrow here goes from left to right} 
 $\phi:X \to Y$, $T^*\phi$ is its {\bf cotangent lift}, defined as the
 canonical relation 
$\{((x,(T\phi)^*(\eta)),(\phi(x),\eta))|x\in X
   \mbox{~and~} \eta \in T^*_{\phi(x)}Y\}$ to $T^*X$ from $T^* Y$.
As a manifold, $T^*f$ may
   be identified with the pulled back vector bundle $\phi^* T^*Y$ over
   $X$.  It is also the image of the conormal bundle to the graph of
   $\phi$ under the symplectomorphism from $T^*(X \times Y)$ to $T^*X
   \times \overline{T^*Y}$ given by reversing the sign of cotangent
     vectors to $Y$.  (A slightly different version of this
 map is called the Schwartz transform
     in \cite{ba-we:lectures}.)  It is easy to check 
 that $T^*$
     embeds $\man$ as a subcategory $T^*\man$ of $\symp$ (in
     particular, that      composition of cotangent lifts is always strongly
     transversal.   

It turns out that we can isolate a property of cotangent lifts which
makes their compositions transversal, and then we can look for more
general situations where this property is satisfied.  Remembering that
the pair $(T^*M,Z_M)$, with $Z_M$ the zero section, is the local model
for any pair consisting of a symplectic manifold and its lagrangian
submanifold, we make the following definition.

\begin{dfn}
\label{dfn-bigmicro}
Let $(X,A)$ and $(Y,B)$ be pairs consisting of a symplectic manifold
and a lagrangian submanifold.   A canonical relation $f$ to $X$ from
$Y$ is {\bf liftlike} with respect to $A$ and $B$ 
if there is a smooth map
$\phi:A\to B$ such that 
\begin{eqnarray*}
f(b) & = & \phi^{-1}(b),\quad\textrm{for all }b\in B,
\\
Tf(v) & = & \big(T\phi\big)^{-1}(v),\quad\textrm{for all }v\in T_{b}B.
\end{eqnarray*}
where $Tf$ is the tangent bundle of $f$ considered as a submanifold of
$TX \times TY$, hence a relation to $TX$ from $TY$.   
\end{dfn}

Every cotangent lift is liftlike with respect to the zero sections,
and the composition of liftlike relations is always transversal near
the relevant lagrangian submanifolds, but in order to get a category,
we must localize around those submanifolds.
 (The use of the prefix ``micro'' below
is meant to correspond to its use by Milnor \cite{mi:microbundles} in
the term ``microbundle''.) 

\begin{dfn}
A {\bf manifold pair}
consists of a manifold $M$ and a closed submanifold $A\subseteq M$. 
Two manifold pairs $(M,A)$ and $(N,B)$ will be considered
equivalent if $A=B$ and if there is a manifold pair
$(U,A)$ such that $U$ is an open subset in both $M$ and $N$ simultaneously.
A {\bf microfold} is an equivalence class $[M,A]$ 
of manifold pairs $(M,A)$.  The (well defined) submanifold $A$ is
the {\bf core} of $[M,A]$.  
\end{dfn}

Note that we require equality of neighborhoods and not
merely diffeomorphism for two manifold pairs to be equivalent.

Most of the standard constructions on manifolds carry over to
microfolds.  In particular, a
{\bf submicrofold} of a microfold $[M,A]$ is a microfold
$[N,B]$ such that $N\subseteq M$ and $B\subseteq A$.
and the product $[M,A] \times [N,B]$ is  $[M \times N, A \times B]$
  A relation between two microfolds is just a submicrofold of their product.
It is (the graph of) a {\bf map} $[M,A] \leftarrow[N,B]$ if it has a
representative which is a map.  This makes the microfolds into a
category $\mic$.  There is a natural forgetful core functor
$[M,A]\mapsto A$ and a cross section thereof $A \mapsto [A,A]$ (with
the obvious actions on morphisms).  

\begin{dfn}
A {\bf symplectic microfold} is a microfold $[M,A]$
together with a germ of symplectic structure on $M$ for which $A$ is
lagrangian.  $[T^*A,A]$ with the canonical symplectic
structure on $T^*A$ is the {\bf cotangent microbundle} of $A$.
A {\bf symplectic micromorphism} to $[M,A]$ from $[N,B]$
is a lagrangian submicrofold of $[M,A] \times [\overline{N},B]$
having a representative which is liftlike with respect to $A$ and
$B$.  The associated map $\phi:A\to B$ is the {\bf core map} of the
micromorphism.  
\end{dfn}

By design, any composition of symplectic micromorphisms is
transversal, so
the symplectic microfolds and micromorphisms form a category
$\micsymp$. 
A symplectic micromorphism is a map if and only it is invertible;
these maps are just the symplectomorphisms in the micro world.
A basic result in \cite{we:symplectic} is (without the microfold
terminology) that every 
symplectic microfold is symplectomorphic to the cotangent microbundle of
its core.  Thus, the restriction of $\mic$ to the cotangent
microbundles is a full subcategory containing the image of the functor
$T^*:\man \to \micsymp$, but with many more morphisms.   There
is also a forgetful functor $\mathbf{CORE} \micsymp \to \man$.

The category $\micsymp$, or rather its extension by a category of
enhanced micromorphisms, carrying half-densities, should be quantized
by a functor to a category of semiclassical 
Fourier integral operators.   What we have so far \cite{ca-dh-we:symplectic2}
is a category 
$\four$ whose objects are manifolds and whose morphisms are certain
operators between smooth half-densities
 on these manifolds which are formal series
in a parameter $\hbar$.  There is a ``wavefront'' functor from $\four$ to 
$\micsymp$ (in fact, to the subcategory whose objects are cotangent
microbundles) for which the inverse image of the identity morphism
over the cotangent microbundle of a manifold $A$ is an algebra of
semiclassical pseudodifferential operators on $A$.  This functor lifts
to a principal symbol functor which attaches a half density to the
wavefront of any operator.   What is missing is a total symbol
calculus which can be make a symbol functor injective.  Also missing
is a construction of operators from general symplectic micromorphisms
which are not acting on cotangent bundles.  Even for cotangent
bundles, the total symbol construction and its inverse appear to
depend on extra structure, such as connections or local coordinates.
The latter allow one to represent symplectic micromorphisms by
generating families which serve as phase functions for the explicit
construction of the kernels of operators as oscillatory integrals.
The general problem seems is somewhat reminiscent of that of passing
from local deformation quantizations of Poisson manifolds
\cite{ko:deformation} to global ones, as in \cite{ca-fe-to:local}.
Perhaps some of the methods of the latter paper will be helpful.

Finally, we note that monoidal objects in the category $\micsymp$
are essentially local symplectic groupoids in the sense of
\cite{we:symplectic groupoids} and correspond to Poisson manifolds.  
The construction of such objects in the formal and analytic categories
was carried out in \cite{ca-dh-fe:formal} and \cite{dh:universal}
using the ``tree-level'' part of the Kontsevich star product.
A good quantization theory for $\micsymp$ should produce algebras from
these monoidal objects.

\section{The Wehrheim-Woodward 2-category}
\label{sec-ww2} 

In the remaining sections, rather than restricting the nature of our canonical
relations, we extend the notion of what a category should be.  We
begin by returning to the linear case.

The second version of Wehrheim and Woodward's quantization follows the
``groupoid philosophy'' 
 that, given an equivalence relation on a set $S$, one should
always try to replace it by the finer structure of a category whose
objects are the elements of $S$ and whose isomorphism classes are the
equivalence classes.   $\symp'$ thus becomes a $2$-category.  In fact,
recalling that 
any morphism $\boldf$ in $\symp'(X,Y)$ is a generalized lagrangian submanifold
in $X\times \overline Y$, we may define the 2-morphism space
$\symp'(\boldf,\boldg)$ to be the Floer cohomology group which
comprises the morphism space to $\boldf$ from $\boldg$ in
$\donsharp(X\times\overline Y)$.   

Similarly, for the category of categories which is the target of $\donsharp$, 
rather than simply identifying natural equivalence classes of
functors, it is appropriate to 
introduce the 2-category structure in which the natural
equivalences become 2-morphisms.  Wehrheim and Woodward now assign to each
Floer cohomology class in $\symp'(\boldf,\boldg)$ a natural transformation
to $\donsharp(\boldf)$ from $\donsharp(\boldg)$ in such a way that
they obtain a 2-functor from the 2-category $\symp'$ of (admissible)
symplectic manifolds, generalized canonical relations, and Floer
cohomology classes, to the 2-category of categories, 
functors, and natural transformations.

The paragraphs above are merely a
schematic description of a {\em tour de force} of symplectic topology
using the authors' theory (see \cite{we-wo:pseudoholomorphic} and
\cite{we-wo:quilted})
of ``quilted pseudoholomorphic curves''.
These are, roughly speaking, piecewise pseudoholomorphic curves
satisfying ``seam'' conditions along smooth (real) curves separating the
smooth pieces of the domain.  But each piece of the curve maps to a
different manifold, with the seams constrained by canonical relations.

\section{Composition of linear canonical relations  as a rational map}
\label{sec-rational}

As was mentioned in the introduction, the composition of linear
canonical relations is not continuous.   For example, if 
$L_1$ and $L_2$
are lagrangian subspaces of $X$,
then $(L_1 \times L_2) \circ (L_2 \times L_1) = L_1 \times L_1$.  If
$L_1$ and $L_2$ are transversal, then 
$L_1 \times L_2$ is the limit
as $a\rightarrow 0$ of the graphs $\Gamma_a = \{(T_ax,x)\}$ of the
symplectomorphisms $T_a$ defined by $T_ax = a^{-1} x$ for $x\in L_1$
and $T_a x = ax$ for $x\in L_2$.   Similarly, $L_2 \times L_1$ is the
limit of $\Gamma _{a^{-1}}$.   The compositions $\Gamma_a \circ
\Gamma_{a^{-1}}$ are all equal to the diagonal $\Gamma_1$ and hence so
is their limit, but the composition of the limits is $L_1\times
L_2$.  

In a paper about spectral
analysis on fractal graphs, Sabot \cite{sa:electrical}introduced, 
in the special case 
of the composition of linear symplectic
reductions with lagrangian subspaces,  a construction which
``fills in'' the discontinuities of the composition operation on
linear canonical relations.
Composition
now becomes multiple valued, just as a
discontinuous step function on the line becomes multiple valued if the
gaps in its graph are filled in with  vertical line segments.  

What
follows below comes from Sabot's construction applied to general
compositions, using the fact (see Section \ref{sec-reduction}) that
composition is also a special case of reduction.  (We save
the details for a future article.)
Although we mostly have the case of vector spaces over $\reals$ in mind,
Sabot works over $\complex$.  In view of the possible more general
applications to finite fields suggested by \cite{gu-ha:geometric}
and \cite{gu-ha:quantization}, we 
will carry out as much as possible of
the construction over an arbitrary field $k$.

For
symplectic vector spaces $X$ and $Y$ over $k$, let $S_{XY}$ denote the
grassmannian (a
manifold when $k$ is $\reals$ or $\complex$)
of all lagrangian subspaces of $X\times Y$ .  For three spaces, $X$,
$Y$, and $Z$ (not necessarily distinct), composition of
linear canonical relations is a mapping $M_{XYZ}:
S_{XY}\times S_{YZ} \to S_{XZ}$.  The restriction of this
mapping to the transversal pairs is continuous, with 
graph
$$T_{XYZ} \defequal \{(f_{XY},f_{YZ},f_{XZ})|f_{XY}\trans f_{YZ} 
\mbox{~and~}
f_{XY}\circ f_{YZ} =f_{XZ}\}\subset S_{XY}\times S_{YZ} \times S_{XZ} .$$ 
$T_{XYZ}$ is dense in the graph of $M_{XYZ}$, but its closure $\tbar_{XYZ}$
in $S_{XY}\times S_{YZ} \times S_{XZ}$ contains more.   Sabot's
description of $\tbar_{XYZ}$ in the case of symplectic reduction
extends to general compositions as follows.   

For for a composable
pair $(f_{XY},f_{YZ})$, we measure their failure to be transversal by
the
 {\bf deficiency}
$d(f_{XY},f_{YZ})$, defined as the codimension of 
$(f_{XY}\times f_{YZ}) \oplus (X \times \Delta_Y \times Z)$ in $X\times
\ybar \times Y \times \zbar$, which is also the dimension of the
intersection $(f_{XY}\times f_{YZ}) \cap (\{0_X\}\times \Delta_Y \times
\{0_Z\})$.  

\begin{thm}
\label{thm-sabotcomposition}
When $k$ is $\reals$ or $\complex$, $\tbar_{XYZ}$ is an algebraic
variety consisting of all triples $(f_{XY},f_{YZ},f_{XZ})$ for
which the codimension of $f_{XZ} \cap f_{XY}\circ f_{YZ}$ in 
$f_{XY}\circ f_{YZ}$ is at most $d(f_{XY},f_{YZ})$.  $T_{XYZ}$ is the
set of its regular points.
\end{thm}

We may think of $\tbar_{XYZ}$ as the graph of a ``continuous multiple
valued function'' $\mbar_{XYZ}$ whose value on $(f_{XY},f_{YZ})$ is a subvariety
$f_{XY} \bullet f_{YZ}$ of $S_{XZ}$ containing the usual composition
$f_{XY} \circ f_{YZ}$.  It is perhaps worth noting that the subvariety $f_{XY}
\bullet f_{YZ}$ is a higher Maslov cycle in the sense of \cite{fu:maslov}.
When $k$ is an arbitrary
field, we may take Theorem \ref{thm-sabotcomposition} as a {\em
  definition} of $\tbar_{XYZ}$ and this multiple valued composition operation

Functions like $\mbar_{XYZ}$ given by relations which are the closure
of graphs of smooth maps are known as ``rational maps'' \cite{gr-ha:principles}.
The next step  (work in progress) is to make
the operation $\bullet$ into the composition operation in a category
$\linsymp$ whose objects are symplectic vector spaces, but which is 
 ``enriched over''
a category $\mathbf{RAT}$ whose objects are algebraic varieties and
whose morphisms are rational maps.  In
other words, the morphism spaces in $\linsymp$ will be the
grassmannians of canonical relations, but the composition operations
will be rational maps.  

\begin{rmk}
\label{rmk-rational}
\emph{
Although the
morphisms in $\mathbf{RAT}$ itself are relations, the composition
operation there is not that of 
$\rel$, but rather the operation which assigns to
rational maps $f$ and $g$ the closure of their set theoretic product.  
To have the composition defined at all, one needs to assume that the
rational maps are ``dominant'' in the sense of having dense range.
This raises the further complication that certain parts of the
structure, in particular the inclusion of the units, cannot be dominant.
It begins to appear that making $\mathbf{RAT}$ into a category
suitable for defining internal categories involves issues similar
the ones we have been dealing with for symplectic categories.
The only way out may be the simplicial approach described in Section
\ref{sec-simplicial}. 
} 
\end{rmk}

As noted in the introduction, quantization of linear
symplectomorphisms is not without its problems.  Although Guillemin and
Sternberg \cite{gu-st:problems} show how to quantize symplectic vector
spaces and to associate operators to canonical relations carrying
half-densities, these operators may be unbounded and cannot always be
composed.   
In fact, the undefined compositions of operators
correspond precisely to nontransversal compositions of canonical
relations.  A simple example is the operator on ``functions'' on
$\reals$ which assigns to each $u$ the product of $u(0)$ with a delta
function at $0$.  This operator, which is a quantization of $L\times
L$ , where $L$ is the
fibre over $0$ in $T^*\reals$, cannot be composed with itself.

The discussion above suggests that we should try to
modify the target category for the usual quantization, either by
introducing a composition of operators which is multiple valued, or by
looking at a category whose morphism spaces, while identified with
spaces of linear operators, themselves admit ``rational maps'' which
may not be defined (or may be multiply-defined) on individual operators.

\begin{rmk}
\emph
{The ``completion'' of the composition of canonical relations to 
something larger also occurs in the microlocal theory of sheaves
\cite{ka-sc:sheaves}.}  
\end{rmk}

\begin{rmk}
\emph
{Segal \cite{se:definition} suggests another approach to building and
quantizing a well-behaved category of linear canonical relations.
For any symplectic vector space $(X,\omega)$ (which could be a product $Y
\times \zbar)$, the positive lagrangian subspaces $L$ of the
complexification $V_\complex$ (with symplectic structure extend from
$V$ by complex bilinearity) are those for which $i\omega(v,\overline
v)>0$ for all nonzero $v\in V$.  Positive canonical relations are
always transversally composable, and their composition is smooth.  They
do not quite form a category, since the identities are missing,
but they may be added ``by hand''.   Quantization of this category
(with the objects enhanced by metaplectic structure and the morphisms
enhanced by half-forms) may be done without any obstructions.  Taking
the real limit remains a problem, though.}
\end{rmk}

\begin{rmk}
\emph
{It would be interesting to see how the linear theory extends
to lagrangian affine subspaces of symplectic affine
spaces.  The hermitian line bundles of geometric quantization may play
a more important role here.}
\end{rmk}

\section{The simplicial picture}
\label{sec-simplicial}

The structure of a category $\calc$ can be encoded in that of a simplicial
object $N(\calc)$ called its {\bf nerve}.  We recall (see for example
\cite{se:classifying}, 
\cite{zh:ngroupoids}, or any book on homological algebra)
 that a simplicial object is a collection of
objects $S^k$ ($k=0,1,2,\ldots$) in some base category (such as sets,
topological 
spaces, or manifolds) together with, 
for each $k$, $k+1$ face morphisms $S^k \to S^{k-1}$
and $k+1$ degeneracy morphisms $S^k\to S^{k+1}$ satisfying the composition
laws of the generators of the category of
order-preserving mappings among  the sets $\Sigma^n =
\{0,1,\ldots,n\}$.   The elements of $S^0$ are sometimes called
{\bf vertices} and those of $S^1$ {\bf edges}.   In $N(\calc)$,
the vertices are the objects of $\calc$ and the edges are the
morphisms, with the faces of a morphism being its target and source,
and the degeneracy operator taking each object to its identity
morphism.  It is 
convenient to write $\calc^{[k]}$ for 
$N^k(\calc)$.  The composition of morphisms is encoded in $\calc^{[2]}$,
whose elements are the composable pairs $(f,g)$, with the face
operators taking each such pair to $f$, $fg$, and $g$, while the
degeneracy operators take $f$ to pairs with an identity morphism
appended to one side or the other.   The rest of the structure is
determined by this part, with $\calc^{[k]}$ being the composable $k$-tuples,
face operators given by the composition of pairs of adjacent entries
or elimination of the entry on one end or the other,
and degeneracy operators by the insertion of identities.  The
associative and identity axioms are equivalent to the compatibility
conditions.    

A simplicial topological (possibly discrete)
space $S$ has a {\bf geometric realization} $|S|$
which is obtained from the disjoint union of $S^k \times \Delta_k$ for
all $k$,
where $\Delta_k$ is the usual $k$-simplex, by gluing them together
using rules derived from the face and degeneracy operators.  The usual
cohomology of $|S|$ is called the cohomology of the simplicial space.
For instance, if $S$ is the nerve $N(G)$ of a group $G$,
$|N(G)|$ is a model for the classifying space $BG$, and its
cohomology is the group cohomology of $G$.  

Not every simplicial object comes from a category or groupoid.  The
ones that do are characterized by so-called ``horn-filling
conditions'', the simplest of which require that 
a pair of edges with a common vertex be fillable (perhaps in a unique
way)  by a 2-simplex of
which they are faces.  
  By weakening these conditions, one arrives at
generalizations of the notions of category and groupoid, as in
\cite{zh:ngroupoids}.  

Let us now apply this idea to the composition
of linear canonical relations.  
For simplicity, we limit our attention to 
linear canonical relations from a fixed
symplectic vector space $X$ to itself.
These form a monoid, i.e. a category $E(X)$ with one object, and we
form the usual nerve in which $E^{[k]}(X)$ is 
the cartesian power $\Lag(X\times \xbar)^{k}$.
 Although each $E^{[k]}(X)$ is a topological
space, and even a smooth manifold, $E^\bullet(X)$ is not a simplicial
topological space because the boundary operators are discontinuous.

The first way to build a simplicial space out of this one is based on
the construction of Wehrheim and Woodward described in Section
\ref{sec-ww1}. 
Within each $E^{[k]}(X)$, there is an open dense subset
$E_\trans^{[k]}(X)$ of
``completely transversal'' sequences, namely
those for which the composition
$f_1\circ\ldots\circ f_k$ is
transversal in the sense that $L_1\times\ldots\times L_k$ is
transversal in $(X\times \xbar)^{k}$ to the multidiagonal
$X\times (\Delta_X)^{k-1}\times \xbar.$  It is not hard to show that this
collection of subsets is invariant under the boundary and degeneracy
operators, and that the restricted operators are smooth, making of
$E^\bullet_\trans(X)$ a simplicial manifold.  

Other constructions might be based on Sabot's multiple-valued
composition.
The simplest one would be to 
take the
closures of the graphs of composition on the 
$E_\trans^{[k]}(X)$ defined in Section \ref{sec-rational} above, but it is not clear that these form
a simplicial object.   It may also be useful to consider structures
 in which there is a whole family of 2-simplices whose edges are a pair $(L_1,L_2)$ of
nontransversely composable relations and an element of $L_1\bullet L_2$.  In any case, it
should be possible to have the 
face and degeneracy operators be ordinary mappings of varieties rather than
rational maps, since the multiple-valuedness of composition is
incorporated in the definition of the spaces of simplices.
Associativity may be taken into account through identification of some
simplices, or through the introduction of extra higher simplices.  

To quantize these simplicial versions of the symplectic category (more
precisely, enhancements thereof), 
one will need to look for similar
simplicial structures derived from the composition of unbounded
operators on function spaces.  This is essentially the point of view
taken in  \cite{ma-we-wo:functoriality}. 
Since
the Floer cohomology
classes which give 2-morphisms in the 2-category described in Section
\ref{sec-ww2} are themselves equivalence classes of cochains, 
Mau, Wehrheim, and Woodward
replace the Donaldson category attached to any admissible symplectic
manifold by a Fukaya-type $A^\infty$ category.  To any admissible
canonical relation they attach an $A^\infty$ functor, and to any Floer
cocycle between such relations a natural transformation of $A^\infty$
functors, in a way which is compatible with composition of strongly
transversal pairs, up to homotopy of $A^\infty$ functors.

\end{document}